\documentclass[12pt,a4paper]{article}
\usepackage{natbib}
\usepackage{amsmath}
\usepackage{amssymb}
\usepackage{amsthm}
\usepackage{graphicx}
\usepackage{microtype}
\usepackage[hidelinks]{hyperref}
\usepackage{verbatim}

\newcommand{\LR}{\mathit{L\mskip-1muR}}
\newcommand{\p}{\theta}
\newcommand{\q}{\bar\theta}
\newcommand{\half}{\tfrac12}
\newcommand{\abs}[1]{\lvert#1\rvert}
\newcommand{\err}{\varepsilon}

\newtheorem{theorem}{Theorem}
\theoremstyle{remark}
\newtheorem*{remark}{Remark}

\author{%
	Wouter Kager\footnotemark[1]%
	\and
	Ronald Meester\thanks{Vrije Universiteit Amsterdam, Department of 
	Mathematics, The Netherlands}}

\title{On the relation between likelihood ratios and $p$-values for testing 
success probabilities of Bernoulli trials}

\begin{document}

\maketitle

\begin{abstract}
	It is well known that there is no direct one-to-one relation between 
	$p$-values and likelihood ratios or Bayes factors, since their relation 
	crucially involves the sample size~$n$. We investigate their (asymptotic) 
	relation in a coin-tossing context where the hypotheses of interest 
	address the success probability of the coin, and where detailed 
	computations are possible. This leads to useful insights in the nature of 
	$p$-values and likelihood ratios. Our results imply, for instance, that 
	under mild conditions, a $p$-value of~$0.05$ cannot correspond to a 
	likelihood ratio larger than~$7.5$, for \emph{any} hypothesis versus a 
	null hypothesis that the success probability has a specific value. We also 
	show it is unlikely one can obtain a large likelihood ratio by tossing a 
	fair coin until the number of heads deviates from the mean by several 
	standard deviations.
\end{abstract}

\medskip\noindent
\textit{Key words and phrases}: $p$-value; likelihood ratio; Bayesian 
analysis; five sigma; evidential value; asymptotics.

\section{Introduction}
\label{s: intro}

It is well known that $p$-values cannot be easily interpreted as a 
quantification of statistical evidence, despite the fact that this is still 
common practice. In \cite{briggs, meesterslootenbook, fluke, royall}, and in 
many other references therein, it is argued that a (very) small $p$-value does 
not at all rule out that the null hypothesis is in fact the best explanation 
of the observed data. Indeed, the fact that the $p$-value is small only means 
that the data (or more extreme outcomes) have a small probability under that 
hypothesis, but this could be the case for any rivaling hypotheses as well. It 
is therefore logically unjustified to dismiss the null hypothesis only because 
the $p$-value is small. Moreover, the mere fact that $p$-values depend on the 
probabilities of unobserved outcomes means they cannot be a proper measure of 
the strength of the evidence against the null hypothesis, as is for instance 
pointed out in~\cite{fluke, royall}.

A likelihood paradigm is often suggested as an alternative (\cite{royall}). In 
that paradigm, one compares the probabilities of the observed data under two 
competing hypotheses. The ratio of these values is called the likelihood 
ratio, or Bayes factor, depending on the circumstances and one's philosophy. 
The likelihood ratio (or Bayes factor) tells us which of the two hypotheses 
explains the data best, and by what factor. As such, it quantifies the 
evidence the data contains in favor of one hypothesis over another, and is a 
genuinely relative notion.

Before we continue and to avoid confusion, a brief discussion of our 
terminology is called for. We use the term \emph{hypothesis} both for a 
`point' hypothesis which fixes the value(s) of the parameter(s), and for the 
situation in which the parameters are thought to follow a certain probability 
distribution. This is consistent with, e.g., \cite{wagenmakers2022} and the 
general use of hypotheses in forensic science, see \cite{meesterslootenbook}. 
For example, in a coin-flip context with unknown success probability~$\p$, we 
will speak of the hypothesis that $\p$ equals~$\frac12$, but also of the 
hypothesis that $\p$ follows a uniform distribution on~$(0,1)$. This is a 
matter of nomenclature only, and our results could be formulated in other 
statistical philosophies as well. 

We use the term \emph{likelihood ratio} not only for the case where all 
parameters in the hypotheses are fixed, but also when these are described by a 
probability distribution and integrated over under either hypothesis. The 
latter is often called a \emph{Bayes factor} in the statistical literature, 
but \cite{edwards1963} for instance use the term \emph{likelihood ratio} for 
both cases as well. Again, this is a matter of convention only. From our point 
of view, the likelihood ratio obtained by an observer will depend on the 
observer's knowledge and understanding of the parameters and is conditional on 
that, rather than on the parameters themselves. The case in which this 
knowledge leads to a point distribution is in this philosophy a special case. 
We refer to \cite{royall} and~\cite{fluke} for gentle introductions to the use 
of likelihood ratios.

Our investigation was initially inspired and motivated by the habit of 
theoretical particle physicists to claim a particle discovery if the data lead 
to a ``five-sigma'' $p$-value (see, e.g., \cite{lyons} and further references 
therein). The ``five sigma'' here refers to a $p$-value that would be obtained 
if an observation of a random variable with a normal distribution were at 
least five standard deviations away from the mean. This boils down to a 
numerical value of roughly~$3\times 10^{-7}$. 

However, even a five-sigma $p$-value does not necessarily imply that the 
relevant likelihood ratios are convincingly large or small. It is very 
important how this five sigma was obtained. To be more precise, it is 
essential to know how many experiments were done, and also why. These warnings 
are well established, of course, but here we look at them from the perspective 
of the asymptotic relation between $p$-values and likelihood ratios. 

We do not work in the extremely complicated world of elementary particle 
physics here, but  in the conveniently controlled environment of coin tosses, 
where explicit asymptotic results can be derived. That is, we consider $n$ 
coin tosses using a coin with unknown success probability~$\p$, and we want to 
make an inference about the value of~$\p$.

In this article, we derive explicit (asymptotic) relations between likelihood 
ratios and $p$-values in a coin-tossing context. It is well known that there 
is no direct one-to-one correspondence between $p$-values and likelihood 
ratios, as their relation involves the sample size~$n$, see, e.g., 
\cite{wagenmakers2022, sellke2001}, and the many references therein. We 
investigate their relation in various situations, varying the hypotheses 
involved and also the experimental design. Some of our results are known (we 
give references in such cases), but are deduced in a new way and in a 
different context. Others (the three theorems and some other results) are new 
as far as we know.

Our analysis also yields some interesting numerical results. For example, we 
show that if $n$ is not too small, a $p$-value of~$0.05$ cannot correspond to 
a likelihood ratio larger than~$7.5$ for \emph{any} hypothesis versus a null 
hypothesis that~$\p$ has a specific value.

The paper is organized as follows. In Section~\ref{s: results} we present our 
main results. We defer the technical proofs of those results to 
Section~\ref{s: proofs}, and we close with a brief discussion of our 
conclusions in Section~\ref{s: discussion}.

\section{Main results}
\label{s: results}

\subsection{Uniform versus fair success probability}
\label{s: H1 vs H2}

We start by comparing the hypothesis~$H_1$ that the unknown success 
probability~$\p$ is uniformly distributed on~$(0,1)$ to the hypothesis~$H_2 : 
\theta=\half$ that the coin is fair. Suppose we toss the coin $n$~times and 
obtain $k$~heads and $n-k$~tails. The likelihood ratio of the data, relative 
to the two hypotheses, is the ratio of the probabilities of the data under 
each hypothesis, namely
\begin{equation}
	\label{LR def}
	\LR
	= \frac{P(\text{data} \mid H_1)} {P(\text{data} \mid H_2)}
	= \frac{\binom{n}{k} \, \int_0^1 \p^k\, (1-\p)^{n-k}\, d\p}
	{\binom{n}{k} \, \left(\frac12\right)^n}.
\end{equation}
Indeed, under $H_1$ the unknown success probability~$\p$ is claimed to be 
uniformly distributed over~$(0,1)$, so we condition on~$\p$ and integrate over 
the density (which is identically equal to~1) to arrive at the expression in 
the numerator. We remark that it is not important whether our data consist of 
the precise sequence of coin tosses or only of the number of heads and tails. 
In the former case, the factor~$\binom{n}{k}$ disappears from both the 
numerator and the denominator, which leads to the same likelihood ratio.

Under $H_2$, the number of heads has mean~$\half n$ with a standard deviation 
of~$\half \sqrt{n}$. We now implicitly define a new parameter~$u$, given $n$ 
and~$k$, via
\[
	k   = \half n - \half u \sqrt{n}, \qquad
	n-k = \half n + \half u \sqrt{n}.
\]
The hereby introduced parameter~$u$ measures how many standard deviations our 
observation is away from the mean under~$H_2$. This will turn out to be a 
useful parametrization. Our interest lies in the behavior of the likelihood 
ratio as a function of $n$ and~$u$, with an emphasis on large~$n$.

In Section~\ref{s: proofs H1 vs H2} we prove the following result: 

\begin{theorem}
	\label{thm: LR bounds}
	Choose $\err \in (0,1)$. If the inequality $\abs{u} \leq \err\sqrt{n}$ 
	holds, then
	\[
		 -\frac12\frac{\err^2}{1-\err^2} - \frac{11}{10n}
		\leq \ln\biggl( \LR\,\sqrt{\frac{2n}{\pi}} \,\,\biggr) -\frac{u^2}{2} 
		- \frac{u^2\err^2}{12} \leq
		\frac{u^2\err^4}{9}.
	\]
\end{theorem}

The takeaway from Theorem~\ref{thm: LR bounds} is that if $n$ is sufficiently 
large given the outcome~$u$ of the experiment, then a good approximation 
for~$\LR$ is
\begin{equation}
	\label{LR approx}
	\LR \approx  \sqrt{ \frac{\pi}{2n} } \, e^{\frac12 u^2}.
\end{equation}
Here, what one considers a ``good'' approximation depends on personal taste 
and is reflected in one's choice of the parameter~$\err$ in the theorem. We 
see from Equation~\eqref{LR approx} that the relation between $\LR$ and~$u$ 
(and therefore between $\LR$ and the $p$-value under~$H_2$, as we explain 
below) crucially involves~$n$. We need both $n$ and~$u$ to compute~$\LR$, and 
hence to assess the relative evidence the data contain for either hypothesis.

We now provide some intuition behind Equation~\eqref{LR approx}. As we show in 
Section~\ref{s: proofs H1 vs H2}, every outcome~$k$ is equally likely under a 
uniformly distributed parameter~$\p$. Therefore, $u$ is roughly uniformly 
distributed over~$(-\sqrt{n},\sqrt{n}\,)$ under~$H_1$, while under~$H_2$, $u$ 
is approximately standard normal. Dividing the uniform density 
over~$(-\sqrt{n},\sqrt{n}\,)$ by~$\varphi(u)$, the standard normal density in 
the point~$u$, reproduces Equation~\eqref{LR approx}. For fixed~$n$, the 
likelihood ratio increases with~$u$ because the outcome becomes increasingly 
unlikely under~$H_2$. Conversely, if we fix~$u$, then the likelihood ratio 
decreases with~$n$ because the outcome becomes increasingly unlikely 
under~$H_1$.

Next we consider $p$-values. The 2-sided $p$-value under~$H_2$ depends only 
on~$u$ (approximately, for large enough~$n$). It is given by
\begin{equation}
	\label{pv def}
	p \approx 2\,\Phi(-\abs{u}\,),
\end{equation}
where $\Phi$ denotes the cumulative distribution function of the standard 
normal distribution. We show in Section~\ref{s: proofs H1 vs H2} that
\begin{equation}
	\label{pv approx}
	p \approx \sqrt{\frac{2}{\pi u^2}} \, e^{-\frac12 u^2},
\end{equation}
and that this leads to the following relation between $\LR$ and~$p$:
\begin{equation}
	\label{LR(p) simple}
	\LR \approx \frac{1}{p\,\sqrt{-2n\ln p}}.
\end{equation}
Clearly, there is no direct one-to-one correspondence between the likelihood 
ratio and the $p$-value, since Equation~\eqref{LR(p) simple} involves the 
sample size~$n$.

\begin{table}
\caption{Likelihood ratios and $p$-values for various choices of $n$ and~$k$, 
	together with the corresponding parameter~$u = (n-\nobreak 2k)/\sqrt{n}$. 
	The table shows the actual values, with our approximations based on 
	Equations \eqref{LR approx} and~\eqref{pv approx} between parentheses.}
\begin{center}
	\begin{tabular}{@{}cccc@{\hskip5pt}cc@{\hskip5pt}c@{}}
		\rule[-6pt]{0pt}{18pt}
		$n$ & $k$ & $u$
		& \multicolumn{2}{c}{$\LR$} & \multicolumn{2}{c}{$p$-value}\\
		\hline
		\rule{0pt}{13pt}
		20 & 6 & 1.789 & 1.288 & (1.388) & 0.1153 & (0.0901)\\
		100 & 40 & 2.000 & 0.913 & (0.926) & 0.0569 & (0.0540)\\
		1000 & 460 & 2.530 & 0.972 & (0.972) & 0.0124 & (0.0129)\\
		10,000 & 4,852 & 2.960 & 1.002 & (1.001)  & 0.00318 & (0.00337)\\
		100,000 & 49,474 & 3.327 & 1.003 & (1.003) & 0.00089 & (0.00095)\\
		1,000,000 & 498,172 & 3.656 & 1.001 & (1.001) & 0.00026 & (0.00027)
	\end{tabular}
\end{center}
\label{t: LR=1}
\end{table}

To illustrate our findings, we have determined, for various values of~$n$, the 
number of heads~$k$ which leads to a likelihood ratio closest to~1. These 
numbers can be found in Table~\ref{t: LR=1}, together with the corresponding 
likelihood ratio and $p$-value. A similar table appears in~\cite{fluke}, and 
\cite{edwards1963} present an analogous table with the $p$-value held fixed 
instead of the likelihood ratio. The results confirm that a $p$-value should 
not be used for evidential purposes: The values of~$k$ in the table constitute 
essentially neutral evidence between $\p$ being uniformly distributed 
on~$(0,1)$ and $\p$ being~$\half$. Nonetheless, the $p$-value under~$H_2$ can 
be very small.

At this point the reader may object that we have been considering only a very 
particular alternative to the hypothesis that $\p$ equals~$\half$. However, we 
shall argue in Section~\ref{s: Hf vs H0} below that similar conclusions can be 
drawn for a wide range of alternative hypotheses.

\begin{figure}
\begin{center}
	\includegraphics{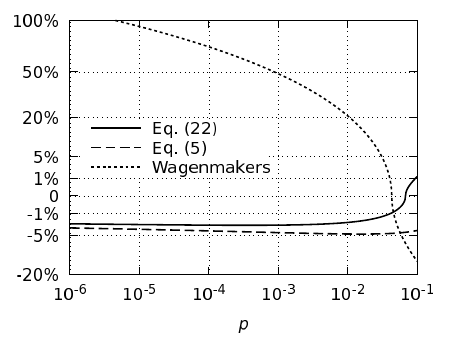}
\end{center}
\caption{Relative errors as a function of the $p$-value of three 
	approximations to the function $(2\sqrt{n}\,\varphi(u))^{-1}$, where $u = 
	\Phi^{-1}(p/2)$: our approximations in Equations \eqref{LR(p) simple} 
	and~\eqref{LR(p)}, and the one based on Wagenmakers' $3p\sqrt{n}$ rule.}
\label{fig: LR(p) rules}
\end{figure}

Equation~\eqref{LR approx} is not new. It also arises in the context of a 
standard normal test statistic, see Equations (2) and~(6) in 
\cite{wagenmakers2022} (where one must take the constant~$A$ in that paper to 
be be~$\sqrt{\pi/2}\,$), and can in this context even be traced back 
to~\cite{jeffreys1936}. Based on this formula, Wagenmakers proposed to use 
$3p\sqrt{n}$ as an approximation to the inverse likelihood ratio~$\LR^{-1}$. 
Equation~\eqref{LR(p) simple} suggests that this ``$3p\sqrt{n}$ rule'' ought 
to be replaced by the rule $p\sqrt{-2n\ln p}$, and an even better 
approximation is derived in Section~\ref{s: proofs H1 vs H2}. Figure~\ref{fig: 
LR(p) rules} compares the three approximations.

\begin{remark}
	Equations \eqref{LR approx}--\eqref{LR(p) simple} are only asymptotically 
	accurate, in the sense that the ratios of the left- and right-hand sides 
	of these equations tend to~1 in the limit~$n\to\infty$, followed in the 
	case of Equations \eqref{pv approx} and~\eqref{LR(p) simple} by the limit 
	$u\to\infty$ or $p\downarrow0$. The reader desiring better approximations 
	for moderate values of $n$ and~$u$ can apply the usual continuity 
	correction for the binomial distribution by setting $v = \abs{u} -\nobreak 
	n^{-1/2}$ and replacing Equations \eqref{pv def}--\eqref{LR(p) simple} by
	\[
		p \approx 2\,\Phi(-v)
		\approx \sqrt{ \frac{2}{\pi v^2} } \, \frac{v^2+1}{v^2+2}
			\, e^{-\frac12 v^2};
	\]
	\[
		\LR \approx \frac{\exp\bigl( \sqrt{-\ln(\pi p^2/2)/n} \,\bigr)}
			{p\sqrt{-n\ln(\pi p^2/2)}}.
	\]
\end{remark}

\subsection{Comparing two point hypotheses}
\label{s: Hx vs H0}

In this section we turn to the case of comparing two given values of the 
unknown success probability. Let $\p_0 \in (0,1)$ be given, and let $H_0$ be 
the hypothesis that $\p$ equals~$\p_0$. Set $\q_0 = 1-\nobreak\p_0$ and $c_0 = 
(\p_0\q_0)^{1/2}$. Suppose we toss the coin $n$~times, and obtain 
$k$~successes. We write $k = \p_0n -\nobreak u\sigma_n$, where $\sigma_n = 
c_0\sqrt{n}$ is the standard deviation corresponding to~$\p_0$. Then we have 
$n -\nobreak k = \q_0n +\nobreak u\sigma_n$, and $u$ measures how many 
standard deviations the outcome is away from the mean if the success 
probability were in fact~$\p_0$.

We would like to compare~$H_0$ with the hypothesis that $\p$ has some other 
specific value. To do this, it turns out to be computationally convenient to 
reparametrize the situation, and consider the likelihood ratio
\[
	\LR_x = \frac{P(\text{data}\mid H_x)}{P(\text{data}\mid H_0)}
\]
relative to the hypotheses $H_0$ and
\[
	H_x : \p = (k-x\sigma_n)/n,
\]
for varying values of~$x$. In Section~\ref{s: proofs Hx vs H0} we prove two 
principal results concerning the likelihood ratio~$\LR_x$, which we formulate 
as Theorems \ref{thm: LRx bounds} and~\ref{thm: LRx max} below.

In Theorems \ref{thm: LRx bounds} and~\ref{thm: LRx max}, $m_0$ denotes the 
largest of the two parameters $\p_0$ and~$\q_0$. The fact that $\p_0$ need not 
be~$\half$ introduces a fundamental asymmetry in the problem. To deal with 
this, we state the two theorems under the assumption that $u$ is positive. 
This simplifies the formulation of the results without losing generality, as 
the case of negative~$u$ can be turned into that of a positive~$u$  by 
interchanging the roles of heads and tails.

\begin{theorem}
	\label{thm: LRx bounds}
	Suppose $u$ is positive and choose $\err \in (0,2/5)$. If the inequality 
	$m_0u \leq \err\sigma_n$ holds, then for all~$x$ in the 
	range~$(-\tfrac32u, \tfrac32u)$ we have
	\[
		\LR_x \leq \exp\Bigl(\, \frac{u^2-x^2}{2}
			+ \abs{\q_0-\p_0} \, \frac{25u^2\err}{12m_0}
			+ \frac{1}{12} \, \frac{u^2\err^2}{1-\err}
			\,\Bigr);\kern1em
	\]
%	and
	\[
		\LR_x \geq \exp\Bigl(\, \frac{u^2-x^2}{2}
			- \abs{\q_0-\p_0} \, \frac{25u^2\err}{12m_0}
			- \frac{297}{32} \, \frac{u^2\err^2}{2-5\err}
			\,\Bigr).
	\]
	Moreover, if either $\p_0\leq\q_0$ or else $250\,(\p_0-\q_0)\,\err \leq 
	12m_0$ holds as well, then we have $\LR_x \leq \exp(-u^2/2)$ for all~$x$ 
	outside the range~$(-\tfrac32u, \tfrac32u)$.
\end{theorem}

This theorem tells us that if $n$ is sufficiently large given the values of 
$u$ and~$\p_0$, then for all practical purposes we may use the following 
approximation for~$\LR_x$ as a function of~$x$:
\begin{equation}
	\label{LRx approx}
	\LR_x \approx \exp\bigl(\, \half u^2 - \half x^2 \,\bigr).
\end{equation}
In particular, the theorem gives us conditions under which we can use this 
approximation to integrate~$\LR_x$ with respect to a probability distribution 
for~$x$, which is what we will do in the next section. These conditions are 
relatively mild when $\p_0$ is close to~$\half$, but can require a very large 
sample size for extreme values of~$\p_0$. It may appear at first sight that 
Equation~\eqref{LRx approx} itself does not contain $n$ and~$\p_0$, but this 
is due to the parametrization we used. Indeed, $x$ measures the deviation from 
the rescaled outcome~$k/n$ in terms of~$\sigma_n$ (which depends on $\p_0$ 
and~$n$) and $n$~itself.

Expression~\eqref{LRx approx} is rather elegant, and illustrates the 
usefulness of the chosen parametrization. It makes mathematically precise a 
phenomenon that was mentioned already in~\cite{edwards1963}, namely that the 
likelihood ratio as a function of~$\p$ becomes sharply peaked around the 
observation of the experiment. Indeed, via $\p = (k-\nobreak x\sigma_n)/n$, we 
see that in terms of~$\p$, the graph of the likelihood ratio is shaped like a 
normal distribution with mean~$k/n$ and standard deviation~$\sigma_n/n = c_0\, 
/ \sqrt{n}$. The likelihood ratio is of course maximal for $x=0$ (or 
equivalently, $\p = k/n$), because this provides the best explanation of the 
data, per construction.

Equation~\eqref{LRx approx} also seems to suggest that $\LR_x$ is bounded 
above by~$\exp(u^2/2)$ for all values of~$x$. It was shown in 
\cite{edwards1963} that this is indeed a uniform upper bound on the likelihood 
ratio in the case of a normally distributed observation, but in the binomial 
case we are dealing with here this is only approximately true. Our next 
theorem provides a uniform upper bound in our context which is sharper than 
the one in Theorem~\ref{thm: LRx bounds}:

\begin{theorem}
	\label{thm: LRx max}
	Suppose $u$ is positive and $\err\in(0,1)$. If $m_0u \leq \err\sigma_n$ 
	holds, then
	\[
		\LR_x \leq \exp\Bigl(\, \frac{u^2}{2}
			+ (\q_0-\p_0)\,\frac{u^2\err}{6m_0}
			+ \frac{u^2}{12}\,\frac{\err^2}{1-\err} \,\Bigr)
		\qquad \text{for all~$x$}.
	\]
\end{theorem}

\begin{table}
\caption{Minimum values of the sample size~$n$ for which the likelihood 
	ratio~$\LR_x$ is uniformly bounded by~7.5 when $u=1.960$, as a function 
	of~$\p_0$.}
\begin{center}
	\begin{tabular}{@{}cccccccccc@{}}
		\rule[-6pt]{0pt}{18pt}%
		$\p_0$ & 0.1 & 0.2 & 0.3 & 0.4 & 0.5 & 0.6 & 0.7 & 0.8 & 0.9 \\
		\hline
		\rule{0pt}{13pt}%
		$n_{\min}$ & 1,527 & 519 & 206 & 74 & 23 & 20 & 24 & 36 & 74
	\end{tabular}
\end{center}
\label{t: nmin}
\end{table}

Theorem~\ref{thm: LRx max} implies that the likelihood ratio for \emph{any} 
alternative hypothesis versus~$H_0 : \p=\p_0$ cannot be much bigger 
than~$\exp(u^2/2)$ if the sample size~$n$ is not too small given $u$ 
and~$\p_0$. This tells us something remarkable in relation to $p$-values: In a 
test of the null hypothesis $\p=\half$, a (two-sided) $p$-value of~0.05 
corresponds to $u=1.960$ and $\exp(u^2/2) \approx 6.8$. Theorem~\ref{thm: LRx 
max} says in this case that \emph{no} hypothesis explains the data better 
than~$H_0$ by much more than a factor~6.8 if $n$ is not too small. How big $n$ 
should be to draw this conclusion depends crucially on the parameter~$\p_0$, 
but the condition in Theorem~\ref{thm: LRx max} is quite mild. To illustrate 
this, Table~\ref{t: nmin} lists the minimum values of~$n$ for which the upper 
bound in the theorem is at most~7.5, for $\p_0$ ranging from 0.1 to~0.9 
(Section~\ref{s: proofs Hx vs H0} explains how these values were computed).

For a likelihood ratio, a value of at most~7.5 is not convincing at all. This 
should make researchers very conservative about using $p$-values in evidential 
context. In fact, in realistic cases the likelihood ratio for an alternative 
hypothesis compared to~$H_0$ will often be (at least) several factors smaller 
than~$\exp(u^2/2)$. Equation~\eqref{LR approx} in Section~\ref{s: H1 vs H2} 
already illustrated this, and we generalize it in the next section. In view of 
the values of~$\exp(u^2/2)$ presented in Table~\ref{t: maxLR}, we conclude 
that also a $p$-value of 0.01, 0.005, or even~0.001, cannot be taken on faith 
as decisive evidence against~$H_0$, as the alternatives might not offer a much 
better explanation of the data than~$H_0$.

\begin{table}
\caption{Values of $u$ and~$e^{u^2/2}$ corresponding to several $p$-values.}
\begin{center}
	\begin{tabular}{@{}ccc@{}}
		\rule[-6pt]{0pt}{18pt}%
		$p$ (one-sided) & $u$ & $e^{u^2/2}$  \\
		\hline
		\rule{0pt}{13pt}%
		0.050 & 1.645 & 3.9 \\
		0.010 & 2.326 & 15.0 \\
		0.005 & 2.576 & 27.6 \\
		0.001 & 3.090 & 118
	\end{tabular}
	\qquad
	\begin{tabular}{@{}ccc@{}}
		\rule[-6pt]{0pt}{18pt}%
		$p$ (two-sided) & $u$ & $e^{u^2/2}$ \\
		\hline
		\rule{0pt}{13pt}%
		0.050 & 1.960 & 6.8 \\
		0.010 & 2.576 & 27.6 \\
		0.005 & 2.807 & 51.4 \\
		0.001 & 3.291 & 224
	\end{tabular}
\end{center}
\label{t: maxLR}
\end{table}

\subsection{Generic likelihood ratio bounds}
\label{s: Hf vs H0}

The likelihood ratio for any alternative hypothesis versus~$H_0$ can in 
principle be computed by integrating the ratio~$\LR_x$ with respect to a 
distribution for~$x$. As we have explained in the previous section, 
Theorem~\ref{thm: LRx bounds} gives us conditions under which we may use the 
approximation in Equation~\eqref{LRx approx} in such computations. These 
conditions boil down to making sure $n$ is sufficiently large given the values 
of $u$ and~$\p_0$, and throughout this section we work under the assumption 
that this is the case and hence that Equation~\eqref{LRx approx} can be used. 

In practice, the alternative hypothesis will typically be formulated in terms 
of a distribution for the parameter~$\p$ rather than for~$x$. So suppose $H_f$ 
is the hypothesis that $\p$ has density~$f(\p)$. We show in Section~\ref{s: 
proofs Hf vs H0} that via $\p = (k -\nobreak x\sigma_n)/n$ and $\sigma_n = 
c_0\sqrt{n}$, and using Equation~\eqref{LRx approx}, it then follows that the 
likelihood ratio for $H_f$ versus~$H_0$ can be computed via
\begin{equation}
	\label{LRf approx}
	\LR_f
%	= \frac{P(\text{data}\mid H_f)}{P(\text{data}\mid H_0)}
	\approx \frac{c_0}{\sqrt n} \, e^{\frac12 u^2}
		\int f\Bigl(\, \frac{k}{n} - x\,\frac{c_0}{\sqrt n} \,\Bigr)
			\, e^{-\frac12 x^2} \, dx.
\end{equation}
For the uniform distribution on~$(0,1)$, for example, this gives
\begin{equation}
	\label{LR uniform}
	\LR_{f=1_{(0,1)}}
	\approx \frac{c_0}{\sqrt n} \, e^{\frac12 u^2}
		\int_{-\infty}^\infty e^{-\frac12 x^2} \, dx
		= c_0\,\sqrt{\frac{2\pi}{n}} \, e^{\frac12 u^2},
\end{equation}
which for $\p_0=c_0=\half$ reproduces Equation~\eqref{LR approx}, as it 
should.

We will now study generic bounds on the likelihood ratio that follow from 
Equation~\eqref{LRf approx} by restricting our attention to certain classes of 
distributions~$f$. To begin, we consider those densities~$f$ under which the 
probability that the distance between $\p$ and the outcome~$k/n$ is at most 
$c_0\abs{u}\, / \sqrt{n}$, is bounded above by~$C\,\abs{u}\, / \sqrt{n}$ for 
some constant~$C>0$. Note that this includes all densities which are bounded 
above by the constant~$C/(2c_0)$ in the mentioned interval around~$k/n$. We 
show in Section~\ref{s: proofs Hf vs H0} that for all such densities we have
\begin{equation}
	\label{LRf bounds}
	\LR_f \lesssim
	1 + C\,\frac{\abs{u}}{\sqrt n} \, e^{\frac12 u^2}.
\end{equation}
This confirms that for a wide range of alternative hypotheses, the likelihood 
ratio is in general much smaller than the value~$\exp(u^2/2)$ shown in 
Table~\ref{t: maxLR}, as we claimed in the previous section. It once again 
points out the danger of interpreting the $p$-value as a measure of evidence 
against~$H_0$.

Next we consider the family of alternative hypotheses under which $\p$ follows 
a normal distribution with mean~$\p_0$ and a fixed standard deviation~$s$. 
Here, we assume $s$ is allowed to vary over a reasonable range given the value 
of~$\p_0$, such as the range~$(0,1-\nobreak\half m_0)$. As we show in 
Section~\ref{s: proofs Hf vs H0}, for the normal density~$f$ with 
parameters~$(\p_0,s^2)$, Equation~\eqref{LRf approx} yields
\begin{equation}
	\label{LR normal}
	\LR_s \approx \frac1{\sqrt{1+\tau_n^2}} \,
		\exp\Bigl(\, \frac{u^2}{2} \, \frac{\tau_n^2}{1+\tau_n^2} \,\Bigr),
\end{equation}
where $\tau_n^2 = ns^2/c_0^2$. One thing we can now do is maximize the 
right-hand side of Equation~\eqref{LR normal} over~$s$. This reproduces what 
we call the \emph{Edwards bound}~$\mathit{EB}$, first derived in 
\cite{edwards1963} for a normal observation, namely
\begin{equation}
	\label{Edwards}
	\sup_s \LR_s
	\approx \mathit{EB}
	= \frac1{\sqrt{eu^2}} \, e^{\frac12 u^2}.
\end{equation}

Using the same argument which lead from Equation~\eqref{LR approx} to 
Equation~\eqref{LR(p) simple}, we show in Section~\ref{s: proofs Hf vs H0} 
that in terms of the $p$-value, the Edwards bound satisfies
\begin{equation}
	\label{EB(p)}
	\mathit{EB} \approx -\frac1{p \ln(\pi p^2/2) \sqrt{\pi e/2}}
\end{equation}
for small~$p$. A similar generic upper bound on the likelihood ratio is the 
bound $(-ep\ln p)^{-1}$ proposed by \cite{sellke2001} and \cite{benjamin2019}, 
and considered first by~\cite{vovk1993}. Equation~\eqref{EB(p)} shows that the 
two bounds are of similar magnitude, the Edwards bound being the smaller of 
the two. In fact, our analysis can be used to prove that the ratio of the two 
bounds converges to the value $\sqrt{e/(2\pi)} \approx 0.66$ in the limit 
$p\downarrow0$.

Researchers should treat these upper bounds with the same wariness as 
$p$-values. Indeed, the likelihood ratio bounds and the $p$-value carry the 
same information, as the relation between them is bijective and does not 
involve the sample size~$n$. Hence, large values of these bounds must not be 
mistaken for a sign of strong evidence against the null hypothesis.

\begin{table}
\caption{The number of coin flips~$n$ and heads~$k$ obtained by two 
	researchers A and~B, together with the value of the parameter~$u$, the 
	$p$-value, the likelihood ratio~$\LR$, and the Edwards 
	bound~$\mathit{EB}$.}
\begin{center}
	\begin{tabular}{@{}c@{\qquad}cc@{}}
		& A & B \\
		\hline
		\rule{0pt}{13pt}%
		$n$ & 10,000 & 100,000,000,000 \\
		$k$ & \phantom04,815 & \phantom049,999,214,176 \\
		$u$ & 3.70 & 4.97
	\end{tabular}
	\qquad
	\begin{tabular}{@{}c@{\qquad}cc@{}}
		& A & B \\
		\hline
		\rule{0pt}{13pt}%
		$p$ & $2.2 \times 10^{-4}$ & $6.7 \times 10^{-7}$ \\
		$\LR$ & 11.8 & 0.916 \\
		$\mathit{EB}$ & 154 & 2,897
	\end{tabular}
\end{center}
\label{t: 5sigma}
\end{table}

An example makes this clear. We consider two researchers A and~B who each toss 
a fair coin a certain number of times~$n$, observing $k$ heads. The numbers 
are given in Table~\ref{t: 5sigma}. In this table, $\LR$ is the likelihood 
ratio from Section~\ref{s: H1 vs H2} for $\p$ uniformly distributed on~$(0,1)$ 
versus $\p=\half$. (A similar table appears in \cite{lyons}.) As it happens, B 
has evidence slightly favoring the null hypothesis, despite his five sigma. We 
also see in the table that the Edwards bound is much larger than the 
likelihood ratio~$\LR$ for both researchers, and is even three orders of 
magnitude larger for researcher~B.

The reason this happens is the maximization procedure used to derive the 
Edwards bound, which is (strongly) biased against the null hypothesis. Indeed, 
our computations show that the standard deviation which maximizes the 
right-hand side of Equation~\eqref{LR normal} is $s = c_0\sqrt{u^2-1}\, / 
\sqrt{n}$. This is (very) small when $n$ is (very) large, making the normal 
distribution which maximizes~$\LR_s$ an atypical member of the considered 
family of alternatives. So although as a family these alternatives may seem 
reasonable, the maximizing member might not be deemed reasonable at all. 
Researchers should be very wary of this when they encounter likelihood ratio 
bounds of this type.

We conclude this section with a brief digression to an amusing example leading 
to a Bayesian interpretation of the $p$-value. Fix $\p_0 \in (0,1)$ as before, 
and choose~$\alpha$ in~$(0,1-\nobreak m_0)$. Let $H^+_\alpha$ and~$H^-_\alpha$ 
be the respective hypotheses that $\p$ is uniform in~$(\p_0, \p_0 +\nobreak 
\alpha)$ or in~$(\p_0 -\nobreak \alpha, \p_0)$. We show in Section~\ref{s: 
proofs Hf vs H0} that
\begin{equation}
	\label{LRalpha}
	\LR_\alpha
	= \frac{P(\text{data}\mid H^+_\alpha)}{P(\text{data}\mid H^-_\alpha)}
	\approx \frac{\Phi(-u)}{1-\Phi(-u)}.
\end{equation}
Here, for $u>0$, we recognize $\Phi(-u)$ as the 1-sided $p$-value for testing 
the null hypothesis~$H_0$. So we have found a one-to-one relation, which does 
not involve~$n$, between the $p$-value and a likelihood ratio, after all. But 
this particular likelihood ratio does not consider the null hypothesis at all, 
and provides a rather unexpected Bayesian interpretation of the $p$-value.

\subsection{Likelihood ratio upon optional stopping}
\label{s: stopping}

In this section we investigate how likely it is that by optional stopping we 
can obtain a small $p$-value, or a large likelihood ratio against the null 
hypothesis that the coin is fair, if the coin is in fact fair. That is, in 
order to obtain a small $p$-value, we want to toss the coin repeatedly until 
the number of heads deviates from the mean by, say, 3 to~5 standard 
deviations. The first question we ask is: How likely is it that after a finite 
number of tosses~$n$ the number of heads~$k$ will satisfy $\abs{k -\nobreak 
\half n} \geq \half c \, \sqrt{n}$, for some given $c\geq1$ measuring the 
deviation from the mean in standard deviations?

We must be careful here. It is well known that runs of consecutive heads or 
tails will occur. As a consequence, we may actually get lucky and achieve the 
desired deviation already after a small number of coin flips. To avoid such 
``beginner's luck'', we assume that the coin is always tossed at least a 
decent number of times (several hundred, say). We call that number~$m$. The 
question then becomes how likely it is we achieve a ``$c$-sigma'' deviation 
(i.e., $c$ standard deviations from the mean) within a finite number of tosses 
\emph{after} first having tossed the coin $m$~times.

This is not a difficult question. The classical Law of the Iterated Logarithm 
tells us that we are \emph{certain} to obtain a $c$-sigma deviation (and hence 
a $p$-value as small as we like) at some point, for any positive~$c$. However, 
the expected number of coin flips required is infinite for any~$c\geq1$, as 
was shown several decades ago; see~\cite{Breiman1967, Shepp1967}. This raises 
the question what we can say about likelihood ratios upon reaching a $c$-sigma 
deviation.

To answer this question, we consider the likelihood ratio from Section~\ref{s: 
H1 vs H2} addressing the hypothesis that $\p$ is distributed uniformly 
on~$(0,1)$ and the null hypothesis that the coin is fair. The value of this 
likelihood ratio at the moment we first reach $c$~standard deviations from the 
mean is a random quantity, because it depends on the random number of coin 
tosses required. In Section~\ref{s: proofs stopping} we argue that the 
expectation of this quantity is given approximately by the surprisingly simple 
formula
\begin{equation}
	\label{E(LR) stopping}
	\frac1{\sqrt{m}} \, \frac{1+c^2}{c}.
\end{equation}

We see that the expected likelihood ratio, at the moment we achieve a 
$c$-sigma deviation, is quite small. This should not come as a surprise given 
that we worked under the assumption that the coin is fair. It is reassuring, 
however, that even this optional stopping procedure is unlikely to produce a 
large likelihood ratio against the null hypothesis, even though it is 
guaranteed to achieve a $c$-sigma deviation, and hence a very small $p$-value.

\section{Technical proofs}
\label{s: proofs}

In this section we have collected detailed proofs and derivations of the 
results in Section~\ref{s: results}. The reader who wishes to skip these 
technical details on a first reading may consider moving on to the discussion 
in Section~\ref{s: discussion}.

In our proofs we make use of the bounds
\begin{equation}
	\label{ln1}
	\frac{x}{1+x} \leq \ln(1+x) \leq x \qquad \text{if $x>-1$}.
\end{equation}
The upper bound here is standard, and the lower bound can be obtained from it 
with the substitution $x\mapsto-x/(1+\nobreak x)$.

\subsection{Proofs for Section~\ref{s: H1 vs H2}}
\label{s: proofs H1 vs H2}

We start by deriving some preliminary bounds. \cite{robbins1955} has proven 
the following version of Stirling's formula, valid for all~$n\geq1$:
\[
	\sqrt{n} \, \exp\Bigl( \frac1{12n+1} \Bigr)
	\leq \frac{n!}{\sqrt{2\pi}} \, \frac{e^n}{n^n} \leq
	\sqrt{n} \, \exp\Bigl( \frac1{12n} \Bigr).
\]
Since by Equation~\eqref{ln1} we also have
\[
	\frac12 \ln\Bigl( 1+\frac1{5n} \Bigr)
	\geq \frac12 \frac{1}{5n+1}
	\geq \frac1{12n} \quad \text{if $n\geq1$},
\]
it follows that for all~$n$ (including $n=0$ if we put $0^0=1$), $n!$ 
satisfies
\begin{equation}
	\label{Stirling}
	\sqrt{\vphantom/n}
	\leq \frac{n!}{\sqrt{2\pi}} \, \frac{e^n}{n^n} \leq
	\sqrt{n+1/5}.
\end{equation}

Next, consider the function
\[
	f(z) = \begin{cases}
		(1-z)\ln(1-z) + (1+z)\ln(1+z) & z\in(-1,1); \\
		2\ln2 & z = -1, 1.
	\end{cases}
\]
From the MacLaurin series of~$\ln(1\pm z)$ we obtain
\[
	f(z)
	= 2\sum_{k=1}^\infty \Bigl(\, \frac1{2k-1} - \frac1{2k} \,\Bigr) z^{2k}
	= 2\sum_{k=1}^\infty \frac{z^{2k}}{2k\,(2k-1)}.
\]
Hence, using Abel's theorem, we conclude that
\begin{align}
	\label{f(z) lower}
	f(z) &\geq z^2 + \tfrac16 z^4; \\
	\label{f(z) upper}
	f(z) &\leq z^2 + \tfrac16 z^4 + \bigl( 2\ln2-1-\tfrac16 \bigr) z^6
		  \leq z^2 + \tfrac16 z^4 + \tfrac29 z^6.
\end{align}

We are now ready to move on to the proof proper of Theorem~\ref{thm: LR bounds}. 
Using the fact that the beta function~$B(a,b)$ satisfies
\[
	B(a,b)
	= \int_0^1 x^{a-1} \, (1-x)^{b-1} \, dx
	= \frac{\Gamma(a) \, \Gamma(b)}{\Gamma(a+b)},
\]
we see that under~$H_1$ the probability of $k$~heads is given by
\[
	P(\text{$k$ heads} \mid H_1)
	= \binom{n}{k} \int_0^1 \p^k \, (1-\p)^{n-k} \, d\p
	= \frac1{n+1}.
\]
Substituting this result into Equation~\eqref{LR def} shows that
\begin{equation}
	\label{LRkn}
	\LR = \frac{k! \, (n-k)!}{n!} \, \frac{2^n}{n+1}.
\end{equation}
Applying the bounds from Equation~\eqref{Stirling} to the factorials in 
Equation~\eqref{LRkn}, we obtain upper and lower bounds on~$\LR$ which are of 
the form
\[
	\frac{\sqrt{2\pi}}{n+1} \, \sqrt{\frac{(k+a)(n-k+a)}{n+b}} \,
	\frac{2^n}{n^n} \, k^k\, (n-k)^{n-k}.
\]
To be more precise, this expression is an upper bound on~$\LR$ if we fill in 
$a=1/5$, $b=0$, and a lower bound if we take $a=0$, $b=1/5$.

We now set $z := u/\sqrt n$, so that
\[
	k   = \half n\, (1-z), \qquad
	n-k = \half n \, (1+z).
\]
This allows us to rewrite our bounds in terms of~$z$ and the function~$f(z)$ 
as
\[
	\sqrt{\frac{\pi\vphantom{()}}{2n\vphantom{()}}} \,
	\sqrt{\frac{(1+2a/n)^2 - z^2}{(1+1/n)^2 \, (1+b/n)}} \,
	\exp\bigl( \half n f(z) \bigr).
\]
For $a=1/5$ and $b=0$, the middle square-root term is bounded above by~1, 
hence by Equation~\eqref{f(z) upper} we obtain the upper bound
\[
	\LR \leq \sqrt{\frac{\pi}{2n}} \,
		\exp\bigl( \half n z^2 \bigl( 1 + \tfrac16 z^2 + \tfrac29 z^4 \bigr) 
		\bigr).
\]
For the lower bound we substitute $a=0$, $b=1/5$ and obtain
\[\begin{split}
	\LR
	&\geq \sqrt{\frac{\pi}{2n}} \,
		\exp\Bigl( \half n f(z) + \frac12\ln(1-z^2)
			- \ln\Bigl( 1+\frac1n \Bigr)
			- \frac12\ln\Bigl( 1+\frac1{5n} \Bigr)
		\Bigr) \\
	&\geq \sqrt{\frac{\pi}{2n}} \,
		\exp\Bigl( \half n z^2 \bigl( 1 + \tfrac16 z^2 \bigr)
			- \frac12 \frac{z^2}{1-z^2} - \frac{11}{10n}
		\Bigr),
\end{split}\]
where the second step uses Equations \eqref{ln1} and~\eqref{f(z) lower}. 
Theorem~\ref{thm: LR bounds} follows.

\begin{figure}
\begin{center}
	\includegraphics{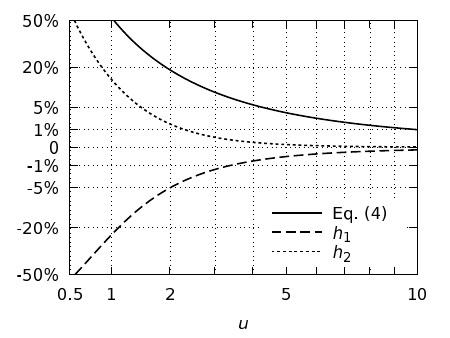}
\end{center}
\caption{The relative errors of three approximations to the 
	function~$2\,\Phi(-u)$: the upper and lower bounds $h_2(u)/A$ 
	and~$h_1(u)/A$, where $A=\sqrt{\pi/2}$, and the approximation in 
	Equation~\eqref{pv approx}.}
\label{fig: errors}
\end{figure}

The approximation for $2\,\Phi(-\abs{u})$ in Equation~\eqref{pv approx} 
follows from the fact that for all positive~$u$ we have $h_1(u) < 
\int_{u}^\infty e^{-x^2/2} \, dx < h_2(u)$, where
\[
	h_1(u) = \frac{u^2}{u^2+1} \, \frac{1}{u} \, e^{-\frac12 u^2}; \qquad
	h_2(u) = \frac{u^2+2}{u^2+3} \, \frac{1}{u} \, e^{-\frac12 u^2}.
\]
This can be shown by comparing the derivatives of the three functions and 
noting that all three functions tend to zero as $u\to\infty$. It can be 
verified numerically that the relative error we make in approximating the 
2-sided $p$-value using Equation~\eqref{pv approx} is less than 19\% when $u$ 
is at least~2, and less than 10\% when $u$ is at least~3; this is illustrated 
in Figure~\ref{fig: errors}.

Squaring Equation~\eqref{pv approx} yields
\[
	u^2\,e^{u^2} \approx \frac{2}{\pi p^2}.
\]
In fact, we know from the previous paragraph that the left-hand side here is 
smaller than the right-hand side. We would like to solve the above equality 
for~$u^2$. To this end, we recall that the principal branch of the Lambert-$W$ 
function is by definition the real-valued function~$W_0(x)$ satisfying
\[
	W_0(x)\,e^{W_0(x)} = x, \qquad x\geq-1/e.
\]
It is not very difficult to show by inspection that for $x>e$ we have
\begin{equation}
	\label{W approx}
	\ln x - \ln\ln x < W_0(x) < \ln x - \ln(\ln x - \ln\ln x).
\end{equation}
To see this, one just has to compute $w(x)\,e^{w(x)}$, taking for~$w(x)$ the 
function on either side of Equation~\eqref{W approx}, and compare the result 
with~$x$.

It follows that in terms of the $p$-value we have approximately
\begin{equation}
	\label{u2(p)}
	u^2 \approx \ln \frac{2}{\pi p^2} - \ln\ln \frac{2}{\pi p^2},
\end{equation}
provided that $p$ is less than $\sqrt{2/(\pi e)} = 0.4839\cdots$. Here we 
chose the left-hand side of Equation~\eqref{W approx} and not the right-hand 
side, because this introduces an error opposite in direction to the one we 
made in Equation~\eqref{pv approx}. Substituting the above approximation 
for~$u^2$ into Equation~\eqref{LR approx} yields
\begin{equation}
	\label{LR(p)}
	\LR \approx \frac{1}{p \, \sqrt{-n\ln(\pi p^2/2)}}.
\end{equation}
When $p$ is small, this may be simplified further to Equation~\eqref{LR(p) 
simple}.

\subsection{Proofs for Section~\ref{s: Hx vs H0}}
\label{s: proofs Hx vs H0}

Recall that $H_0$ is the hypothesis $\{\p=\p_0\}$, and $H_x$ is the hypothesis 
that the success probability~$\p$ equals $(k-\nobreak x\sigma_n)/n$, where $k 
= \p_0n -\nobreak u\sigma_n$ and $\sigma_n^2 = \p_0\q_0 n$. We assume 
throughout this section that~$u$ is positive, but the results we obtain can be 
translated to the case $u<0$ by symmetry considerations.

We write $u +\nobreak x = y = z \sigma_n$. Then under~$H_x$ we have $\p = \p_0 
-\nobreak \p_0\q_0 z$, hence the likelihood ratio for $H_x$ versus~$H_0$ 
becomes
\[
	\LR_x
	= \frac{P(\text{data} \mid H_x)}{P(\text{data} \mid H_0)}
	= \Bigl( \frac{\p_0-\p_0\q_0z}{\p_0} \Bigr)^k \,
		\Bigl( \frac{\q_0+\p_0\q_0z}{\q_0} \Bigr)^{n-k}.
\]
We take the logarithm of this expression, and regard it as a function~$g(y)$ 
of the variable~$y$. Substituting $k = \p_0n -\nobreak u\sigma_n$, this 
function can be written in terms of the variable~$z = y/\sigma_n$ as
\begin{equation}
	\label{lnLRx}
	g(y) = (\p_0n-u\sigma_n) \ln(1-\q_0z) + (\q_0n+u\sigma_n) \ln(1+\p_0z),
\end{equation}
where the range is determined by the inequalities $-\q_0<\p_0\q_0z<\p_0$. By 
the nature of the problem, the function~$g(y)$ should be maximal for $y=u$, 
which corresponds to $x=0$, and should be decreasing in~$y$ for $y>u$, and 
increasing for~$y<u$. This is confirmed by the derivative of~$g(y)$, which is 
given by
\begin{equation}
	\label{g'(y)}
	g'(y)
	= \frac{\q_0u-\sigma_n}{1-\q_0z} + \frac{\p_0u+\sigma_n}{1+\p_0z}
	= \frac{u-y}{(1-\q_0z)(1+\p_0z)}.
\end{equation}

Applying Equation~\eqref{ln1} to~\eqref{lnLRx} yields the simple upper bound
\begin{equation}
	\label{g(y) bound}
	g(y)
	\leq -\q_0z\,(\p_0n-u\sigma_n) + \p_0z\,(\q_0n+u\sigma_n)
	= uy.
\end{equation}
More accurate bounds on~$g(y)$ can be found by comparing $g(y)$ with the 
function
\begin{equation}
	\label{h(y)}
	h(y) = uy-\half y^2 + (\q_0-\p_0) \, \Bigl(
		\frac{uy^2}{2\sigma_n} - \frac{y^3}{3\sigma_n} \Bigr),
\end{equation}
which has derivative
\begin{equation}
	\label{h'(y)}
	h'(y) = (u-y) + (\q_0-\p_0)\,(u-y)\,z.
\end{equation}
Combining Equations \eqref{g'(y)} and~\eqref{h'(y)} yields
\begin{align}
	g'(y)-h'(y)
	&= (u-y) \, \frac{1 - (1-\q_0z)(1+\p_0z)(1+(\q_0-\p_0)\,z)}
		{(1-\q_0z)(1+\p_0z)} \notag \\
	&= (u-y) \, z^2 \, \frac{(\q_0-\p_0)^2+\p_0\q_0+\p_0\q_0\,(\q_0-\p_0)\,z}
		{(1-\q_0z)(1+\p_0z)}.
	\label{g'-h'}
\end{align}
Now consider the fraction in Equation~\eqref{g'-h'}. From $-\q_0 < \p_0\q_0z < 
\p_0$ one sees that its numerator take values between $\p_0^2$ and~$\q_0^2$. 
Hence, if by $m_0$ we denote the largest of $\p_0$ and~$\q_0$, then the 
numerator is positive and at most~$m_0^2$, while the denominator is at least 
$1-\nobreak m_0\abs{z}$.

It follows that for all~$y$ in the range~$(-u,u)$ we have
\[
	0\leq g'(y)-h'(y) \leq
		(u-y) \, \frac{y^2}{\sigma_n^2} \, \frac{m_0^2}{1-m_0u/\sigma_n}.
\]
Integration of the right-hand side over~$y$ from $-u$ to~$0$ now yields
\begin{equation}
	\label{bounds y<0}
	0\geq g(y)-h(y) \geq -\frac{7}{12} \, \frac{m_0^2}{1-m_0u/\sigma_n} \, 
	\frac{u^4}{\sigma_n^2},
	\qquad -u\leq y\leq 0,
\end{equation}
while integration over the range~$(0,u)$ gives
\begin{equation}
	\label{bounds 0<y<u}
	0\leq g(y)-h(y) \leq \frac1{12} \, \frac{m_0^2}{1-m_0u/\sigma_n} \, 
	\frac{u^4}{\sigma_n^2},
	\qquad 0\leq y\leq u.
\end{equation}
Finally, we consider the range~$(u,cu)$ for $c>1$. Observe that in this range 
the derivative of $g(y)-h(y)$ is negative, so the previous upper bound still 
applies. Furthermore, Equation~\eqref{g'-h'} in this case tells us that
\[
	g'(y)-h'(y) \geq
		(u-y) \, \frac{y^2}{\sigma_n^2} \, \frac{m_0^2}{1-cm_0u/\sigma_n},
\]
and integration over the range~$(u,cu)$ yields
\begin{equation}
	\label{bounds u<y}
	g(y)-h(y) \geq -\frac{1+3c^4-4c^3}{12} \, \frac{m_0^2}{1-cm_0u/\sigma_n} 
	\, \frac{u^4}{\sigma_n^2},
	\qquad u\leq y\leq cu.
\end{equation}

Since~$g(y)$ has an absolute maximum in~$y=u$, Theorem~\ref{thm: LRx max} now 
follows directly from the upper bound in Equation~\eqref{bounds 0<y<u} and the 
fact that
\[
	h(u) = \frac{u^2}{2} + (\q_0-\p_0) \, \frac{u^3}{6\sigma_n}.
\]
It also follows that the upper bound on~$g(u)$ equals~$\ln7.5$ if the 
equality
\[
	\frac{u^2}{2} + (\q_0-\p_0) \, \frac{u^3}{6\sigma_n} + \frac1{12} \, 
	\frac{m_0^2u^4}{\sigma_n^2-m_0u\sigma_n} = \ln7.5
\]
holds. Multiplying this expression by $\sigma_n^2 - m_0u\sigma_n$ yields a 
quadratic equation in~$\sigma_n$. From the larger of the two solutions of this 
equation (for $u=1.960$) we obtain the values of~$n_{\min}$ in Table~\ref{t: 
nmin}.

We move on to Theorem~\ref{thm: LRx bounds}. Note that since $y=u+x$, the 
range~$(-\tfrac32u, \tfrac32u)$ for~$x$ corresponds to the range~$(-\half u, 
\tfrac52u)$ for~$y$. (We chose this range based on symmetry considerations, 
but other choices are possible here.) From the definition of~$h(y)$ in 
Equation~\eqref{h(y)} we get
\begin{equation}
	\label{h(5u/2)}
	h\bigl( \tfrac52u \bigr)
	= - \frac{5u^2}{8} - (\q_0-\p_0) \, \frac{25u^3}{12\sigma_n},
\end{equation}
and it is easily seen, based on Equation~\eqref{h'(y)}, that the term $25u^3 / 
(12\sigma_n)$ on the right is an upper bound on the absolute value of the 
function
\[
	\frac{uy^2}{2\sigma_n} - \frac{y^3}{3\sigma_n}
\]
over the specified range for~$y$. Since furthermore the expression $uy-\half 
y^2$ in Equation~\eqref{h(y)} is equal to $\half (u^2-x^2)$ in terms of~$x$, 
the first two bounds on~$\LR_x$ in Theorem~\ref{thm: LRx bounds} now follow 
from Equations~\eqref{bounds y<0}--\eqref{bounds u<y}, taking $c=5/2$.

We now consider $y$ outside the specified range. Equation~\eqref{g(y) bound} 
already gives $g(y) \leq -\half u^2$ for $y \leq -\half u$, so it only remains 
to consider $y \geq \tfrac52u$. We work under the assumption that 
$m_0u/\sigma_n \leq \err < 2/5$ holds. Then the upper bound in 
Equation~\eqref{bounds 0<y<u} is at most~$u^2/45$, hence
\[
	g(y)
	\leq g\bigl( \tfrac52u \bigr)
	\leq h\bigl( \tfrac52u \bigr) + \frac{u^2}{45}.
\]
From Equation~\eqref{h(5u/2)} it is clear that the right-hand side is smaller 
than~$-\half u^2$ when $\q_0$ is at least~$\p_0$. Furthermore, for $\q_0$ 
smaller than~$\p_0$ we have
\[
	g(y) \leq - \frac{5u^2}{8} + (\p_0-\q_0) \, \frac{25u^2\err}{12m_0}
		+ \frac{u^2}{45}.
\]
If the second term on the right is at most~$u^2/10$, then the whole right-hand 
side is again smaller than~$-\half u^2$. This completes the proof of 
Theorem~\ref{thm: LRx bounds}.

\subsection{Proofs for Section~\ref{s: Hf vs H0}}
\label{s: proofs Hf vs H0}

Let $H_f$ be the hypothesis represented by the density~$f(\p)$, and recall 
that $\p = (k -\nobreak x\sigma_n)/n$ and $k = \p_0n -\nobreak u\sigma_n$ in 
our chosen parametrization. Then
\[
	f(\p)\,d\p
	= -\frac{\sigma_n}{n} \,
		f\Bigl(\, \frac{k}{n} - x \, \frac{\sigma_n}{n} \Bigr)\, dx,
\]
hence using $\sigma_n = c_0\sqrt{n}$, Equation~\eqref{LRf approx} follows from 
Equation~\eqref{LRx approx} by integration with respect to the density~$f$. If 
$f$ is the indicator function of the interval~$(0,1)$, then 
Equation~\eqref{LRf approx} gives
\[
	\LR_f
	\approx \frac{c_0}{\sqrt n} \, u^{\frac12 u^2}
		\int_{-u-\q_0\sqrt{n}/c_0}^{-u+\p_0\sqrt{n}/c_0} e^{-\frac12 x^2}
		\, dx,
\]
from which Equation~\eqref{LR uniform} follows if $n$ is large enough.

Now suppose the density~$f$ is such that under~$H_f$, the probability 
that~$\p$ deviates at most $\abs{u} \, \sigma_n/n$ from~$k/n$ is bounded above 
by~$C\,\abs{u} / \sqrt{n}$. Then
\[
	\int_{\abs{x}>\abs{u}} f(\p) \, e^{\frac12 u^2-\frac12 x^2} \, d\p
	\leq \int_{\abs{x}>\abs{u}} f(\p) \, d\p
	\leq 1,
\]
while
\[
	\int_{\abs{x}\leq\abs{u}} f(\p) \, e^{\frac12 u^2-\frac12 x^2} \, d\p
	\leq e^{\frac12 u^2} \int_{\abs{x}\leq\abs{u}} f(\p) \, d\p
	\leq C\, \frac{\abs{u}}{\sqrt n} \, e^{\frac12 u^2}.
\]
Adding these inequalities yields Equation~\eqref{LRf bounds}.

Next, suppose $f$ is the normal density with mean~$\p_0$ and standard 
deviation~$s$. Then under~$H_f$, $x$ has a normal distribution with mean~$-u$ 
and standard deviation~$\tau_n = ns/\sigma_n$. We therefore have
\[
	\LR_s
	\approx \frac{1}{\sqrt{2\pi\tau_n^2}} \,
		\int \exp\biggl(\, \frac{u^2}{2} -\frac{x^2}{2}
			- \frac{(x+u)^2}{2\tau_n^2}
		\,\biggr) \, dx.
\]
Completing the square yields
\[
	\LR_s
	\approx \frac{1}{\sqrt{2\pi\tau_n^2}} \,
		\int \exp\biggl(\, \frac{u^2}{2} \, \frac{\tau_n^2}{1+\tau_n^2}
			-\frac{1+\tau_n^2}{2\tau_n^2} \,
				\Bigl( x + \frac{u}{1+\tau_n^2} \Bigr)^2
			\,\biggr) \, dx,
\]
and evaluating this integral produces Equation~\eqref{LR normal}.

By taking the derivative of the right-hand side of Equation~\eqref{LR normal} 
with respect to~$\tau_n^2$, it is not difficult to verify that $\LR_s$ is 
maximal when $\tau_n^2$ is equal to~$u^2-1$. In terms of~$s$, this means the 
maximum is attained when $s^2$ equals $c_0^2\,(u^2-\nobreak1)/n$. Substituting 
$\tau_n^2 = u^2-1$ into Equation~\eqref{LR normal} yields the Edwards bound, 
Equation~\eqref{Edwards}. Using Equation~\eqref{u2(p)} then produces
\[
	\mathit{EB}
	\approx \frac1{p\sqrt{-(\pi e/2) \ln(\pi p^2/2)}} \,
		\frac1{\sqrt{\,\ln(\pi p^2/2) - \ln\ln(\pi p^2/2)}}.
\]
For small~$p$ we may omit the log-log term, which yields 
Equation~\eqref{EB(p)}.

Finally, let $\alpha \in (0,1 -\nobreak m_0)$, and let $H^+_\alpha$ be the 
hypothesis that $\p$ is uniform on $(\p_0, \p_0 +\nobreak \alpha)$. Note that 
this corresponds to $x$ being uniform between $-u - \alpha n/\sigma_n$ 
and~$-u$. Therefore, for large enough~$n$,
\[
	\frac{P(\text{data} \mid H^+_\alpha)}{P(\text{data} \mid H_0)}
	\approx \frac{\sigma_n}{\alpha n} \, e^{\frac12 u^2}
		\int_{-\infty}^{-u} e^{-\frac12 x^2} \, dx
	= \frac{\sigma_n \sqrt{2\pi}}{\alpha n} \, e^{\frac12 u^2} \, \Phi(-u).
\]
Similarly, for the hypothesis $H^-_\alpha$ that $\p$ is uniform 
on~$(\p_0-\nobreak\alpha, \p_0)$ we obtain
\[
	\frac{P(\text{data} \mid H_{-\alpha})}{P(\text{data} \mid H_0)}
	\approx \frac{\sigma_n}{\alpha n} \, e^{\frac12 u^2}
		\int_{-u}^\infty e^{-\frac12 x^2} \, dx
	= \frac{\sigma_n \sqrt{2\pi}}{\alpha n} \, e^{\frac12 u^2}
		\, \bigl( 1-\Phi(-u) \bigr).
\]
Equation~\eqref{LRalpha} follows by taking the ratio of these two likelihood 
ratios.

\subsection{Proofs for Section~\ref{s: stopping}}
\label{s: proofs stopping}

Let $N$ denote the random number of coin tosses (including the first, fixed 
number of $m$~tosses) until the number of heads first deviates by at least~$c$ 
standard deviations from its mean. If we stop tossing at that point, then 
according to Equation~\eqref{LR approx} we have obtained a likelihood ratio 
given by
\[
	\LR_N \approx \sqrt{\frac{\pi}{2N}} \, e^{\frac12 c^2}.
\]
This is a random variable determined by~$N$, and we are interested in its 
mean, and therefore in computing~$E(N^{-1/2})$.

The discrete nature of the problem makes this computation hard, but we can 
obtain a decent approximation (if $m$ is large) by considering the analogous 
stopping problem for Brownian motion. That is, let $W_t$ be standard Brownian 
motion and set
\[
	T := \inf\bigl\{\, t\geq m : \abs{W_t} \geq c\sqrt{t} \,\bigr\}.
\]
Then $E(T^{-1/2})$ should be a good approximation to~$E(N^{-1/2})$. The 
Brownian motion case was studied in \cite{Shepp1967}. Shepp showed that
\[
	E(\,T^\mu \mid W_m=0)
	= \frac{ m^\mu } { M\bigl(-\mu,\frac12,\frac12 c^2\bigr) },
\]
where $M(a,b,z)$ is Kummer's confluent hypergeometric function.

This result unfortunately involves conditioning on $W_m$ being zero, which in 
our coin-tossing setting means conditioning on the number of heads after the 
initial~$m$ tosses being equal to its mean. However, it turns out that Shepp's 
calculation can be modified to compute the unconditional 
expectation~$E(T^\mu)$ as well. We omit the details, but we claim that by 
following the steps of Shepp it can be shown that
\begin{equation}
	\label{Shepp}
	E(T^\mu)
	= m^\mu \, \biggl(\, 2\, \Phi(-c)
		+ \sqrt{\frac2\pi} \,
			\frac{ M\bigl(1-\mu, \frac32, \frac12 c^2\bigr) }
				{ M\bigl(-\mu, \frac12, \frac12 c^2\bigr) }
			\, c \, e^{-\frac12 c^2} \,\biggr).
\end{equation}

For positive~$\mu$ the function $M(-\mu,\half,z)$ has a positive root, which 
means that~$E(T^\mu)$ is finite only if~$c$ is sufficiently small. To be 
precise, it follows from tabulated values of the smallest positive root 
of~$M(a,b,z)$ that $E(T)$ is finite only if $c<1$, while $E(T^{1/2})$ is 
finite only if $c<1.3069$ \cite[Table~13.2] {abramowitzstegun}. The same 
conclusions then hold for $E(N)$ and~$E(N^{1/2})$ if $T$ is comparable to the 
number of coin tosses~$N$ in our discrete setting. Justification for this is 
given in \cite{Breiman1967}, where it is shown that in the conditional case 
mentioned above, the random variables $T$ and~$N$ do have the same tail 
behavior.

We conclude that for $c$ larger than~1.3069, the inverse likelihood 
ratio~$\LR_N^{-1}$ has infinite mean. As for the mean of~$\LR_N$ itself, we 
note that in the case $\mu = -\half$, the confluent hypergeometric functions 
in Equation~\eqref{Shepp} simplify considerably \cite[Section~13.6] 
{abramowitzstegun}, which leads to
\[
	E(T^{-1/2})
	= \frac{1}{\sqrt{m}} \, \Bigl(\, 2\,\Phi(-c)
		+ \sqrt{2/\pi} \, c \, e^{-\frac12 c^2} \,\Bigr).
\]
If we now assume $E(N^{-1/2}) \approx E(T^{-1/2})$ and approximate the term 
$2\,\Phi(-c)$ using Equation~\eqref{pv approx}, we obtain the expression in 
Equation~\eqref{E(LR) stopping}.

\section{Discussion and conclusion}
\label{s: discussion}

This contribution reinforces the well-known thesis that reporting $p$-values 
without context is statistically rather dangerous. It has been observed, for 
instance in~\cite{ioannidis}, that using $p$-values in an evidential role 
easily gives rise to wrong results and claims. In the current article, we have 
analyzed the relation between $p$-values and likelihood ratios in detail, and 
we have concluded that they display very different asymptotic behavior in 
terms of the observation~$u$ and sample size~$n$. Our results not only confirm 
that this is the case, but also shed some light on the underlying reasons.

Of course, our results were only obtained in the controlled environment of 
coin flips, and they do not verbatim apply to other circumstances. However, we 
note two things: First, the point we want to make is also philosophical in 
nature. Although in different circumstances it may no longer be reasonable to 
expect that exact expressions such as ours can be obtained, the fact that 
$p$-values and likelihood ratios scale differently will probably still hold; 
there is no reason why they would not. Second, especially in applied sciences 
like sociology, psychology, and medicine, the context is quite often one of 
repeated experiments for which the coin-tossing context is entirely 
appropriate. 

There are therefore both philosophical and practical reasons to be very 
cautious with the use of $p$-values. In particular, the fact that in the 
coin-tossing context a $p$-value of~0.05 cannot correspond to a likelihood 
ratio larger than~7.5 under mild conditions on the sample size should make all 
researchers very suspicious about this methodology. The same caution should be 
observed regarding the likelihood ratio (or Bayes factor) bounds that have 
been proposed in the literature. As we have demonstrated, the connection 
between the strength of the evidence as measured by a likelihood ratio and a 
$p$-value is non-trivial and rather complicated.

\bibliographystyle{apalike}
\bibliography{References}

\begin{thebibliography}{}

\bibitem[Abramowitz and Stegun, 1965]{abramowitzstegun}
Abramowitz, M. and Stegun, I.~A., editors (1965).
\newblock {\em Handbook of Mathematical Functions with Formulas, Graphs, and
  Mathematical Tables}.
\newblock Dover Publications, Inc., New York.

\bibitem[Benjamin and Berger, 2019]{benjamin2019}
Benjamin, D. and Berger, J. (2019).
\newblock Three recommendations for improving the use of $p$-values.
\newblock {\em The American Statistician}, 73:186--191.

\bibitem[Breiman, 1967]{Breiman1967}
Breiman, L. (1967).
\newblock First exit times from a square root boundary.
\newblock In {\em Proc. {F}ifth {B}erkeley {S}ympos. {M}ath. {S}tatist. and
  {P}robability ({B}erkeley, {C}alif., 1965/66), {V}ol. {II}: {C}ontributions
  to {P}robability {T}heory, {P}art 2}, pages 9--16. Univ. California Press,
  Berkeley, CA.

\bibitem[Briggs, 2016]{briggs}
Briggs, W. (2016).
\newblock {\em Uncertainty, the {S}oul of {M}odeling, {P}robability and
  {S}tatistics}.
\newblock Springer.

\bibitem[Edwards et~al., 1963]{edwards1963}
Edwards, W., Lindman, H., and Savage, L. (1963).
\newblock Bayesian statistical inference for psychological research.
\newblock {\em Psychological Review}, 70:193--242.

\bibitem[Ioannidis, 2005]{ioannidis}
Ioannidis, J. (2005).
\newblock Why most publised research findings are false.
\newblock {\em Plo{S} {M}edicine}, 2:e124.

\bibitem[Jeffreys, 1936]{jeffreys1936}
Jeffreys, H. (1936).
\newblock Further significance tests.
\newblock {\em Mathematical {P}roceedings of the {C}ambridge {P}hilosophical
  {S}ociety}, 32:416--445.

\bibitem[Lyons, 2013]{lyons}
Lyons, L. (2013).
\newblock Discovering the significance of 5 sigma.
\newblock arXiv:1310.1284 (preprint).

\bibitem[Meester and Slooten, 2020]{meesterslootenbook}
Meester, R. and Slooten, K. (2020).
\newblock {\em Theory and Philosophy of Statistical Evidence in Forensic
  Science}.
\newblock Cambridge University Press.

\bibitem[Meester and Slooten, 2021]{fluke}
Meester, R. and Slooten, K. (2021).
\newblock {\em Fact or fluke}.
\newblock Amsterdam {U}niversity {P}ress.

\bibitem[Robbins, 1955]{robbins1955}
Robbins, H. (1955).
\newblock A remark on {S}tirling's formula.
\newblock {\em The American Mathematical Monthly}, 62:26--29.

\bibitem[Royall, 1996]{royall}
Royall, R. (1996).
\newblock {\em Statistical Evidence}.
\newblock CRC Press.

\bibitem[Sellke et~al., 2001]{sellke2001}
Sellke, T., Bayarri, M., and Berger, J. (2001).
\newblock Calibration of $p$ values for testing precise null hypotheses.
\newblock {\em The American Statistician}, 55:62--71.

\bibitem[Shepp, 1967]{Shepp1967}
Shepp, L.~A. (1967).
\newblock A first passage problem for the {W}iener process.
\newblock {\em Annals of Mathematical Statistics}, 38:1912--1914.

\bibitem[Vovk, 1993]{vovk1993}
Vovk, V. (1993).
\newblock A logic of probability, with applications to the foundations of
  statistics.
\newblock {\em Journal of the Royal Statistical Society, Series~B},
  55:317--351.

\bibitem[Wagenmakers, 2022]{wagenmakers2022}
Wagenmakers, E. (2022).
\newblock Approximate objective {B}ayes {F}actors from $p$-values and sample
  size: the $3p\sqrt{n}$ rule.
\newblock {\em Psy ArXiv preprints (\@osf.io)}.

\end{thebibliography}

\end{document}